\documentclass[8pt,a4paper]{article}
\usepackage{amsfonts}
\usepackage{amsmath}
\usepackage{amssymb}
\usepackage{amsthm}
\usepackage[dvips]{graphicx}
\usepackage{latexsym}
\usepackage{mathrsfs}
\usepackage{amscd}
\usepackage{indentfirst}
\usepackage{textcomp}
\usepackage[numbers,square,sort&compress]{natbib}
\usepackage[colorlinks=true,allcolors=black ]{hyperref}
\allowdisplaybreaks[4]

\newtheorem{Theorem}{Theorem}[section]

\newtheorem{Corollary}[Theorem]{Corollary}
\newtheorem{Lemma}[Theorem]{Lemma}

\theoremstyle{definition}
\newtheorem{Definition}{Definition}[section]
\numberwithin{equation}{section}

\DeclareMathOperator{\Der}{Der}
\DeclareMathOperator{\BDer}{BDer}
\DeclareMathOperator{\IBDer}{IBDer}

\DeclareMathOperator{\Span}{span}
\renewcommand{\d}{\mathrm{d}}
\renewcommand{\div}{\mathrm{div}}
\newcommand{\BF}{\mathbb{F}}

\newcommand{\BN}{\mathbb{N}}
\newcommand{\BZ}{\mathbb{Z}}

\begin{document}
\title{{\bf  Skew-symmetric super-biderivations of the special odd Hamiltonian superalgebra $SHO(n,n;\underline{t})$}}
\author{\normalsize \bf Da  Xu$^1$\,\,\,\,Xiaoning  Xu$^1$}
\date{{{\small{   School of Mathematics and Statistics, Liaoning University, Shenyang, 110036,
China  }}}}
\maketitle

\begin{abstract}
This paper aims to study the skew-symmetric super-biderivations of the special odd Hamiltonian superalgebra $SHO(n,n;\underline{t})$. Let $HO$ denote the odd Hamiltonian Lie superalgebra $HO(n,n;\underline{t})$ and $SHO$  the special odd Hamiltonian Lie superalgebra $SHO(n,n;\underline{t})$ over a field of characteristic $p > 2$. Utilizing the abelian subalgebra $T_{SHO}$ and the weight space decomposition of $HO$ with respect to $T_{SHO}$, we show the action of a skew-symmetric super-biderivation on the elements of $T_{SHO}$ and some specific elements of $SHO$. Moreover, we prove that all skew-symmetric super-biderivations of $SHO$ are inner.
\end{abstract}
\renewcommand{\thefootnote}{\fnsymbol{footnote}}
\footnote[0]{ Project supported by National Natural Science
Foundation of China (No.11501274) and the general project for the
educational department of Liaoning Province, China(No. LJKMZ20220454, LJKZ0096)..
\\ *Corresponding author: xuxiaoning@lnu.edu.cn (X. Xu)}

\section{Introduction}

The theory of Lie superalgebras over a field of characteristic zero has seen a remarkable progress, both in mathematics and pyhsics \cite{a9,a10,a16,a17,a18}.
The complete classification of the finite-dimensional simple modular Lie superalgebras remains an open problem \cite{a11}. The research of modular Lie superalgebra has undergone an evolution during the last ten years, especially in the classical simple modular Lie superalgebras and the structures and representations of simple modular Lie superalgebra of Cartan type. The eight families of finite-dimensional Cartan-type modular Lie superalgebras $W, S, H, K, HO, KO, SHO$ and $SKO$ were discussed in \cite{a5,a12,a13,a14,a25}. The superderivation algebras, second cohomologies, filtrations, representations of the eight families of finite-dimensional simple modular Lie superalgebras of Cartan type were investigated (see \cite{a7,a15,a25}, for example).

As is well known, the study of derivations is very active because of their importance in Lie algebras and Lie superalgebras. With further research about the theory of derivations, it is therefore natural to begin the investigations of biderivations and commuting maps on Lie algebras \cite{a1,a2,a4,a8,a19,a20,a21}.
In particular, the notations of super-biderivations and skew-symmetric super-biderivations was introduced in \cite{a6, a22}. The results on super-biderivations and skew-symmetric super-biderivations of classical simple Lie superalgebras  arose in \cite{a23}.
The skew-symmetric super-biderivations  of generalized Witt Lie superalgebra $W(m,n;\underline{t})$ were proved to be inner in \cite{a3}. In \cite{a24}, there were similar results for Contact Lie superalgebra $K(m,n;\underline{t})$.

This paper is devoted to studying the skew-symmetric super-biderivations of the special odd Hamiltonian superalgebra $SHO(n,n;\underline{t})$. And the paper is arranged as follows. In Section 2, we review the basic definitions concerning the special odd Hamiltonian superalgebra $SHO(n,n;\underline{t})$. In Section 3, we give the definition of skew-symmetric super-biderivations of Lie superalgebras. And some useful conclusions about the super-biderivations are given. In Section 4, we use the method of the weight space decomposition of $HO(n,n;\underline{t})$ with respect to the abelian subalgebra $T_{SHO}$ to prove that all skew-symmetric super-biderivations of $SHO(m,n;\underline{t})$ are inner.

\section{Preliminaries}

In this section, a brief review of the special odd Hamiltonian Lie superalgebras $SHO(n,n;\underline{t})$ is given (see \cite{a12}).

Hereafter $\BF$ is an algebraically closed field of characteristic $p>2$ and $\BZ_2 = \{\bar{0}, \bar{1}\}$ is the additive group of order 2. For a vector superspace $V=V_{\bar{0}}\, \oplus\, V_{\bar{1}}$, we write $\d(x)=\alpha$ for the parity of $x\in V_\alpha,\alpha\in \BZ_2$. If $V=\oplus_{i\in \BZ}V_i$ is a $\BZ$-graded vector space, for $x\in V_i,\,i\in \BZ$, $x$ is a $\BZ$-homogeneous element and its $\BZ$-degree $i$ . Throughout this paper, we should mention that once the symbol $\d(x)$  appears in an expression, it implies that $x$ is a $\BZ_2$-homogeneous  element.

Let $\BN_+$ be the set of positive integers and $\BN$ the set of non-negative integers.
Given $n\in\BN_+\backslash\{1\}$. For $\alpha=(\alpha_{1},\alpha_{2},\cdots,\alpha_{n})\in \BN^n$, put $|\alpha|=\sum_{i=1}^{n}\alpha_{i}$. For $\beta=(\beta_{1},\beta_{2},\cdots,\beta_{n})\in \BN^n$, we write $\alpha+\beta=(\alpha_{1}+\beta_{1},\alpha_{2}+\beta_{2},\cdots,\alpha_{n}+\beta_{n})$,  $\binom{\alpha}{\beta}=\prod_{i=1}^n\binom{\alpha_i}{\beta_i}$, $\alpha\leq\beta\Longleftrightarrow\alpha_i\leq\beta_i, i=1,2,\cdots,n$. For $\varepsilon _i:=( \delta _{i1},\delta _{i2},\cdots ,\delta _{in})$, where $\delta_{ij}$ is the Kronecker symbol, we abbreviate $x^{(\varepsilon_{i})}$ to $x_{i}$, $i=1, 2, \cdots, n$.  We call $\mathcal{U}(n)$ a $divided\ power\ algebra$ which denotes the  $\BF$-algebra of power series in the variable $x_{1},x_{2},\cdots,x_{n}$. The following formula holds in $\mathcal{U}(n)$:
  \[x^{(\alpha)}x^{(\beta)}=\binom{\alpha+\beta}{\alpha}x^{(\alpha+\beta)},\ \forall\ \alpha,\beta\in \BN^n.\]

Let $\Lambda(n)$ denote the $Grassmann\ superalgebra$ over $\BF$ in $n$ variables $x_{n+1},x_{n+2},\cdots,x_{2n}$. Denote the tensor product $\mathcal{U}(n)\otimes\Lambda(n)$ by $\Lambda(n,n)$. Then $\Lambda(n,n)$ have  a $\BZ_2$-gradation induced by the trivial $\BZ_2$-gradation of  $\mathcal{U}(n)$ and the natural $\BZ_2$-gradation of $\Lambda(n)$:
  \[\Lambda(n,n)_{\bar{0}}=\mathcal{U}(n)\otimes\Lambda(n)_{\bar{0}},\qquad \Lambda(n,n)_{\bar{1}}=\mathcal{U}(n)\otimes\Lambda(n)_{\bar{1}}.\]

Obviously, $\Lambda(n,n)$ is an associative superalgebra. For $g\in\mathcal{U}(n), f\in\Lambda(n)$, we simply write $g\otimes f$ as $gf$. The following formulas hold in $\Lambda(n,n)$:
  \[x_ix_j=-x_jx_i, \qquad i,j=n+1,\cdots,2n.\]
  \[x^{(\alpha)}x_j=x_jx^{(\alpha)}, \quad\ \forall\ \alpha\in \BN^n,j=n+1,\cdots,2n.\]
For $k = 1,\cdots,n$, set
  \[B_k:=\{\langle i_1,i_2,\cdots,i_k\rangle\mid n+1\leq i_1<i_2<\cdots<i_k\leq 2n\}\]
and $B(n):=\bigcup_{k=0}^nB_k$, where $B_0=\varnothing$. For $u=\langle i_1,i_2,\cdots,i_k\rangle\in B_k$, set $|u|=k, \{u\}=\{i_1,i_2,\cdots i_k\}$ and $x^u=x_{i_1}x_{i_2}\cdots x_{i_k}$. Specially, let $|\varnothing|=0$, $x^\varnothing=1$. It is obvious that $\{x^{(\alpha)}x^u\mid \alpha\in \BN^n,u\in B(n)\}$ is an $\BF$-basis of $\Lambda(n,n)$.

Let $Y_0=\{1,2,\cdots,n\}$, $Y_1=\{n+1,\cdots,{2n}\}$, and $Y=Y_0\cup Y_1$. Let $D_1,D_2,\cdots,D_{2n}$ be the linear transformations of $\Lambda(n,n)$ such that
  \begin{align*}
  	D_i(x^{(\alpha)}x^u)=
  	\begin{cases}
  		x^{(\alpha-\varepsilon_i)}x^u,\quad&\forall\ i\in Y_0,\\
  		x^{(\alpha)}\partial_i(x^u),\quad &\forall\ i\in Y_1,
  	\end{cases}
  \end{align*}
where $\partial_i$ is the special derivation of $\Lambda(n)$. Then $D_1,D_2,\cdots,D_{2n}$ are superderivations of the superalgebra $\Lambda(n,n)$, and $\d(D_i)=\tau(i)$, where
  \begin{align*}
  	\tau(i)=
  	\begin{cases}
  		\bar{0},\qquad \forall\ i\in Y_0,\\
  		\bar{1},\qquad \forall\ i\in Y_1.
  	\end{cases}
  \end{align*}
Let
  \[W(n,n):=\bigg\{\sum_{i=1}^{2n}f_iD_i\mid f_i\in \Lambda(n,n),\ \forall\ i\in Y\bigg\}.\]
Then $W(n,n)$ is an infinite-dimensional Lie superalgebra which is contained in $\Der(\Lambda(n,n))$ and the following formula holds:
  \[ [fD_i,gD_j]=fD_i(g)D_j-(-1)^{\d(fD_i)\d(gD_j)}gD_j(f)D_i, \]
for all $f,g\in\Lambda(n,n)$ and $i,j\in Y$.

Fix two $n$-tuples of positive integers $\underline{t}=(t_1,t_2,\cdots,t_n)$ and $\pi=(\pi_1,\pi_2,\cdots,\pi_n)$,
where $\pi_i=p^{t_i}-1$ for all $i\in Y_0$ and $p$ is the characteristic of the basic field $\BF$. Let
  \[\Lambda(n,n;\underline{t}):= \Span_{\BF}\{x^{(\alpha)}x^u\mid \alpha\in A(n,\underline{t}),u\in B(n)\},\]
where $A(n,\underline{t} ) =\{\alpha=(\alpha _1,\alpha _2,\cdots ,\alpha _n)\in \BN^n\mid 0\leq \alpha _i\leq \pi _i,i\in Y_0 \}$. Then $\Lambda(n,n;\underline{t})$ is a subalgebra of $\Lambda(n,n)$. Let
  \[W(n,n;\underline{t}):=\bigg\{\sum_{i=1}^{s}f_iD_i\mid f_i\in \Lambda(n,n;\underline{t}),\ \forall\ i\in Y\bigg\}.\]
Then $W(n,n;\underline{t})$ is a finite-dimensional simple Lie superalgebra, which is called the generalized Witt Lie superalgebra. Note that $W(n,n;\underline{t})$ possesses a $\BZ$-graded structure:
  \[W(n,n;\underline{t})=\bigoplus\limits_{r=-1}^{\xi-1}W(n,n;\underline{t})_r,\]
by letting $W(n,n;\underline{t})_r:=\Span_{\BF}\{x^{(\alpha)}x^uD_j\mid|\alpha|+|u|=r+1,j\in Y\}$ and $\xi:=|\pi|+n$.

Put
\begin{align*}
	i'&=\begin{cases}
		i+n,&\quad  1\leq i\leq n,\\
		i-n,& \quad n< i\leq 2n. \\
	\end{cases}
\end{align*}
Define a linear mapping $T_{H}: \Lambda(n,n;\underline{t})\rightarrow W(n,n;\underline{t})$  such that
  \[T_{H}(f)=\sum_{i=1}^{2n}(-1)^{\tau(i)\d(f)}(D_i(f)D_{i'}),\]
for all $f\in \Lambda(n,n;\underline{t})$. Note that $T_H$ is odd and that
  \[[ T_H(f),T_H(g)]=T_H(T_H(f)(g)),\]
for all $f,g\in \Lambda(n,n;\underline{t})$. Let
  \[HO(n,n;\underline{t}):=\Span_{\BF}\{T_H(f)\mid f\in \Lambda(n,n;\underline{t})\}.\]
Then $HO(n,n;\underline{t})$ is a finite-dimensional simple Lie superalgebra, which is called the odd Hamiltonian superalgebra. 
We define $S'(n,n;\underline{t})$ by means of the following mapping, called divergence:
\begin{align*}
	\div:&\begin{cases}
		W(n,n;\underline{t})&{\rightarrow} \,\,\,\,\Lambda(n,n;\underline{t}),\\
		\,\,\,\,\,\,\,f_iD_i &{\mapsto}\,\,\,\, (-1)^{\tau(i)\d(f_i)}D_i(f_i). \\
	\end{cases}
\end{align*}
Note that $\div$ is even and that
\[\div([D,E])=[\div(D),E]+[D,\div(E)],\]
for all $D,E\in W(n,n;\underline{t})$. Let
\[S'(n,n;\underline{t}):=\{D\in W(n,n;\underline{t})\mid \div(D)=0\}.\]
Then $S'(n,n;\underline{t})$ is a $\BZ$-graded subalgebra of $W(n,n;\underline{t})$. Let
\begin{align*}
	SHO'(n,n;\underline{t})&:=S'(n,n;\underline{t})\cap HO(n,n;\underline{t}),\\
	\overline{SHO}(n,n;\underline{t})&:=[SHO'(n,n;\underline{t}),SHO'(n,n;\underline{t})],\\
	SHO(n,n;\underline{t})&:=[\overline{SHO}(n,n;\underline{t}),\overline{SHO}(n,n;\underline{t})].
\end{align*}
Then $SHO(n,n;\underline{t})$ is a finite-dimensional simple Lie superalgebra, which is called the special odd Hamiltonian superalgebra. Note that $SHO(n,n;\underline{t})$ possesses a $\BZ$-graded structure:
  \[SHO(n,n;\underline{t})=\bigoplus\limits_{r=-1}^{\xi-5}SHO(n,n;\underline{t})_r,\]
  where $$SHO(n,n;\underline{t})_r:=\Span_{\BF}\{T_{H}( x^{( \alpha)}x^{u})\mid|\alpha|+|u|=r+2,\,\mathrm{div}(T_{H}( x^{( \alpha)}x^{u}))=0\}.$$
 For convenience, we use $HO$, $SHO$ and $SHO_r$ denote $HO(n,n;\underline{t})$, $SHO(n,n;\underline{t})$ and its $\BZ$-graded subspace $SHO(n,n;\underline{t})_r$ respectively.

\section{The notion of super-biderivation}
 
 In this section,  considering the $\BZ_2$-degree of bilinear map, we give the definition of skew-symmetric super-biderivations of Lie superalgebras, which  is different from that mentioned in \cite{a6}. Moreover, we obtain some useful conclusions about the skew-symmetric super-biderivations.

Let $G$ be a Lie algebra over an algebraically closed field. Recall that a linear map $D : G \rightarrow G$ is a derivation of $G$ if
\[D([x,y])=[D(x),y]+[x,D(y)],\]
for all $x, y\in G$. And we call a bilinear map $\varphi : G \times G\rightarrow G$ is a biderivation of $G$ if the following axioms are satisfied:
\begin{align*}
	\varphi (x,[y,z])&=[\varphi (x,y),z]+[y,\varphi (x,z)],\\
	\varphi([x,y],z)&=[\varphi(x,z),y]+[x,\varphi(y,z)],
\end{align*}
for all $x, y, z\in G$. A biderivation $\varphi$ is called skew-symmetric if $\varphi(x, y) = -\varphi(y, x)$ for all $x, y \in G$.

A Lie superalgebra is a vector superspace $L=L_{\bar{0}}\,\oplus\, L_{\bar{1}}$ with an even bilinear mapping $[\cdot,\cdot]:L\times L\rightarrow L$ satisfying the following axioms:
\begin{align*}
	[x,y]&=-(-1)^{\d(x)\d(y)}[y,x],\\
	[x,[y,z]]&=[[x,y],z]+(-1)^{\d(x)\d(y)}[y,[x,z]],
\end{align*}
for all $x,y,z\in L$. Recall that a linear map $D:L\rightarrow L$ is a superderivation of $L$ if
\[D([x,y])=[D(x),y]+(-1)^{\d(D)\d(x)}[x,D(y)],\]
for all $x,y\in L$. Meanwhile, we write $\Der_{\bar{0}}(L)$ (resp. $\Der_{\bar{1}}(L)$) for the set of all superderivations of $\BZ_2$-degree $\bar{0}$ (resp. $\bar{1}$) of $L$. Let $\Der(L)=\Der_{\bar{0}}(L)\oplus\Der_{\bar{1}}(L).$

A $\BZ_2$-homogeneous bilinear map $\phi$ of $\BZ_2$-degree $\gamma$ of $L$ is a bilinear map such that $\phi(L_\alpha,L_\beta)\subset L_{\alpha+\beta+\gamma}$ for any $\alpha,\beta\in\BZ_2$. Specially, we call $\phi$ is even if $\d(\phi)=\gamma=\bar{0}$.

\begin{Definition}\label{def3.1}
	We call a bilinear map $\phi : L\times L\rightarrow L$ is a skew-symmetric super-biderivation of L if the following axioms are satisfied:
	\begin{align}
\phi(x,[y,z])&=[\phi(x,y),z]+(-1)^{(\d(\phi)+\d(x))\d(y)}[y,\phi(x,z)], \label{eq3.1}\\
\phi(x,y)&=-(-1)^{\d(\phi)\d(x)+\d(\phi)\d(y)+\d(x)\d(y)}\phi(y,x),\nonumber
	\end{align}
	for all $x,y,z\in L$. 
\end{Definition}
 Observe that if $\phi$ is a skew-symmetric super-biderivation, the equation below is obviously fulled. 
\begin{align}
&\phi([x,y],z)=[x,\phi(y,z)]+(-1)^{(\d(\phi)+\d(z))\d(y)}[\phi(x,z),y]. \label{eq3.2}
\end{align}
Meanwhile, we write $\BDer_{\bar{0}}(L)$ (resp. $\BDer_{\bar{1}}(L)$) for the set of all skew-symmetric super-biderivations of $\BZ_2$-degree $\bar{0}$ (resp. $\bar{1}$) of $L$. Let $\BDer(L)=\BDer_{\bar{0}}(L)\oplus\BDer_{\bar{1}}(L).$

\begin{Lemma}\label{lem3.1}
	If the bilinear map $\phi_{\lambda}$ : $L\times L\rightarrow L$ is defined by
	\[\phi_{\lambda}( x,y ) =\lambda [ x,y],\]
	for all $x,y\in L$, where $\lambda \in \BF$, then $\phi_{\lambda}$ is a skew-symmetric super-biderivation of $L$. This class of super-biderivations is called inner. Write $\IBDer(L)$ for the set of all inner super-biderivations of $L$.
	
	\begin{proof}
		Due to $\phi_{\lambda}( x,y ) =\lambda [ x,y]$, we deduce that $\phi_{\lambda}$ is an even super-biderivation, i.e. $\d(\phi_\lambda)=\bar{0}$. By the skew-symmetry of Lie superalgebras, it is easy to seen that
		\[\phi_\lambda(x,y)=-(-1)^{\d(\phi_\lambda)\d(x)+\d(\phi_\lambda)\d(y)+\d(x)\d(y)}\phi_\lambda(y,x),\]
		for any $x, y \in L$. And it is easily obtain the following three equalities:
		\begin{align*}
			\phi_\lambda(x,[y,z])&=\lambda[x,[y,z]],\\
			[\phi_\lambda(x,y),z]&=\lambda[[x,y],z],\\
			(-1)^{(\d(\phi_\lambda)+\d(x))\d(y)}[y,\phi_\lambda(x,z)]&=(-1)^{\d(x)\d(y)}\lambda[y,[x,z]].
		\end{align*}
		Since the graded Jacobi identity $[x,[y,z]]=[[x,y],z]+(-1)^{\d(x)\d(y)}[y,[x,z]]$,  we have that
		\[\phi_\lambda(x,[y,z])=[\phi_\lambda(x,y),z]+(-1)^{(\d(\phi_\lambda)+\d(x))\d(y)}[y,\phi_\lambda(x,z)],\]
		for any $x,y,z \in L$. Similarly, it follows that
		\[\phi_\lambda([x,y],z)=[x,\phi_\lambda(y,z)]+(-1)^{(\d(\phi_\lambda)+\d(z))\d(y)}[\phi_\lambda(x,z),y],\]
		for any $x, y, z \in L$. The proof is completed.
	\end{proof}
	
\end{Lemma}

\begin{Lemma}\label{lem3.2}		
	Let $L$ be a Lie superalgebra. Suppose that $\phi$ is a skew-symmetric super-biderivation of $L$. Then
	\[[\phi(x,y),[u,v]]=(-1)^{\d(\phi)(\d(y)+\d(u))}[[ x,y] ,\phi( u,v)],\]
	for any $x,y,u,v\in L$.
	
	\begin{proof}
		Due to the Definition \ref{def3.1}, there are two different ways to compute $\phi([x, u], [y, v])$. From the equation \eqref{eq3.1}, it follows that
		\begin{align*}
			\phi ([x,u],[y,v])
			=\ &[\phi([x,u],y), v]+(-1)^{(\d(\phi)+\d(x)+\d(u))\d(y)}[y,\phi([x,u],v)] \\
			=\ &[[x,\phi(u,y)], v]+(-1)^{(\d(\phi)+\d(y))\d(u)}[[\phi(x, y), u], v] \\
			&+(-1)^{(\d(\phi)+\d(x)+\d(u))\d(y)}[ y,[ x,\phi (u, v)]] \\
			&+(-1)^{(\d(\phi)+\d(x)+\d(u))\d(y)+(\d(\phi)+\d(v))\d(u)}[y, [\phi(x, v), u]].
		\end{align*}
		According to the equation \eqref{eq3.2}, one gets		
		\begin{align*}
			\phi([x,u],[y,v])
			=\ &[x,\phi(u,[y,v])]+(-1)^{(\d(\phi)+\d(y)+\d(v))\d(u)}[\phi(x,[y,v]),u] \\
			=\ &[x,[\phi(u,y),v]]+(-1)^{(\d(\phi)+\d(u))\d(y)}[x,[y,\phi(u,v)]] \\
			&+(-1)^{(\d(\phi)+\d(y)+\d(v))\d(u)}[[\phi(x,y),v],u] \\
			&+(-1)^{(\d(\phi)+\d(y)+\d(v))\d(u)+(\d(\phi)+\d(x))\d(y)}[[y, \phi(x, v)], u].
		\end{align*}
		
		Comparing two sides of the above two equations, and using the graded Jacobi identity of Lie superalgebras, we have that
		\begin{equation}\label{eq3.3}
			\begin{aligned}
				&[\phi(x,y),[u,v]]-(-1)^{\d(\phi)(\d(y)+\d(u))}[[x,y],\phi(u,v)]\\
				=\ &(-1)^{\d(y)\d(u)+\d(y)\d(v)+\d(u)\d(v)}([\phi(x,v),[u,y]]-(-1)^{\d(\phi)(\d(v)+\d(u))}[[x,v],\phi(u, y)]).
			\end{aligned}		
		\end{equation}		
		And we set
		\[f(x,y;u,v)=[\phi(x,y),[u,v]]-(-1)^{\d(\phi)(\d(y)+\d(u))}[[x,y],\phi(u,v)].\]
		According to the equation \eqref{eq3.3}, it easily seen that
		\[f(x,y;u,v)=(-1)^{\d(y)\d(u)+\d(y)\d(v)+\d(u)\d(v)}f(x,v;u,y).\]		
		On the one hand, we have
		\begin{align*}
			f(x,y;u,v)&=-(-1)^{\d(u)\d(v)}f(x,y;v,u)\\
			&=-(-1)^{\d(u)\d(v)}(-1)^{\d(y)\d(v)+\d(y)\d(u)+\d(v)\d(u)}f(x,u;v,y)\\
			&= (-1)^{\d(y)\d(u)}f(x,u;y,v).
		\end{align*}
		On the other hand, we also get
		\begin{align*}
			f(x,y;u,v)&=(-1)^{\d(y)\d(u)+\d(y)\d(v)+\d(u)\d(v)}f(x,v;u,y)\\
			&=-(-1)^{\d(y)\d(u)+\d(y)\d(v)+\d(u)\d(v)}(-1)^{\d(u)\d(y)}f(x,v;y,u)\\
			&=-(-1)^{\d(u)\d(y)}f(x,u;y,v).
		\end{align*}		
		Hence, it follows at once that
		\[f(x,y;u,v)=-f(x,y;u,v).\]		
		Due to the $\mathrm{char}(\BF)\neq 2$, we obtain that $f(x,y;u,v)=0.$ Hence
		\[[\phi(x,y),[u,v]]=(-1)^{\d(\phi)(\d(y)+\d(u))}[[x,y],\phi(u, v)].\]
		The proof is completed.
	\end{proof}
	
\end{Lemma}

\begin{Lemma}\label{lem3.3}		
	Let $L$ be a Lie superalgebra. Suppose that $\phi$ is a skew-symmetric super-biderivation of $L$. If $\d(x)+\d(y)=\bar{0}$, then
	\[[\phi ( x,y ) ,[ x,y]] =0,\]
	for any $x,y\in L$.
	
	\begin{proof}
		By virtue of Lemma \ref{lem3.2}, we can have
		\[[\phi(x,y),[x,y]]=(-1)^{\d(\phi)(\d(y)+\d(x))}[[x,y],\phi(x, y)].\]		
		In view of $\d(x)+\d(y)=\bar{0}$, it is following that
		\begin{align*}
			[\phi(x,y),[x,y]]&=[[x,y],\phi(x,y)]\\
			&=-[\phi(x,y),[x,y]].
		\end{align*}
		Therefore, we obtain $[\phi(x,y),[x, y] ]=0$.
	\end{proof}
	
\end{Lemma}

\begin{Lemma}\label{lem3.4}   	
	Let $L$ be a Lie superalgebra. Suppose that $\phi$ is a skew-symmetric super-biderivation of $L$. If $[x,y]=0$, then
	\[\phi(x,y)\in Z_L([L,L]),\]
	where $Z_L([L,L])$ is the centralizer of  $[L,L]$. 	
	
	\begin{proof}
		If $[x,y]=0$, then we obtain
		\begin{align*}
			[\phi(x,y),[u,v] ]&= (-1)^{\d(\phi)(\d(y)+\d(u))}[[x,y],\phi[u,v]]\\
			&= 0,
		\end{align*}
		for any $u, v \in L$. So we get that $\phi(x,y)\in Z_L([ L,L])$.
	\end{proof}
	
\end{Lemma}

\section{Skew-symmetry super-biderivation of $SHO(n,n;\underline{t})$}

In this section, we use the method of the weight space decomposition of $HO$ with respect to the abelian subalgebra $T_{SHO}$ to prove that all skew-symmetric super-biderivations of $SHO$ are inner. For convenience, we use $A$ and $B$ denote $A(n,\underline{t})$ and $B(n)$.

Set $T_{SHO}=\Span_{\BF}\{T_H(x_ix_{i'}-x_jx_{j'})\mid i,j\in Y_0\}$. It is easy to see that $T_{SHO}$ is an abelian subalgebra of $SHO$. Obviously, $T_{SHO}$ is also an abelian subalgebra of $HO$. For any $T_H(x^{(\alpha)}x^u)\in HO$, we have
  \[[T_{H}( x_ix_{i'}-x_jx_{j'}) ,T_{H}( x^{( \alpha)}x^{u})] = (\delta_{( i'\in u)}- \alpha_i+\alpha_j-\delta_{( j'\in u)}) T_{H}( x^{( \alpha)}x^u),\]
where
  \[\delta_{(\text{P})}=
  \begin{cases}
  	\ 1 \qquad \text{P is ture},\\
  	\ 0 \qquad \text{P is false}.
  \end{cases}\]
Fixed $\alpha\in A$ and $u\in B$, we define a linear function $(\alpha+\langle u\rangle):T_{SHO}\rightarrow \BF$
  \[(\alpha+\langle u\rangle)(T_{H}( x_ix_{i'}-x_jx_{j'}))=\delta_{( i'\in u)}- \alpha_i+\alpha_j-\delta_{( j'\in u)}.\]
Furthermore, $HO$ have a weight space decomposition with respect to $T_{SHO}$:
\[HO=\bigoplus_{(\alpha+\langle u\rangle)}HO_{(\alpha+\langle u\rangle)},\]
where
  \begin{align*}
  	HO_{(\alpha+\langle u\rangle)}=\Span_{\BF}\{T_H(x^{(\beta)}x^v)\in HO \mid &[T_{H}( x_ix_{i'}-x_jx_{j'}),T_H(x^{(\beta)}x^v)]=\\
  	&(\delta_{( i'\in u)}- \alpha_i+\alpha_j-\delta_{( j'\in u)}) T_H(x^{(\beta)}x^v)\}.
  \end{align*}
For any $T_H(x^{(\alpha)}x^u)\in SHO$, let
\[SHO_{(\alpha+\langle u\rangle)}:=HO_{(\alpha+\langle u\rangle)}\cap SHO.\]

\begin{Lemma}\label{lem4.1}  	
	Let $\phi$ be a skew-symmetric super-biderivation of $SHO$. If $[x,y] = 0$ for $x,y \in SHO$, we have
	\[\phi(x,y) = 0.\]	
	
	\begin{proof}
		Since $SHO$ is a simple Lie superalgebra, it is obvious that $SHO=[SHO,SHO]$ and $Z(SHO)=0$. According to Lemma \ref{lem3.4}, if $[x,y]=0$ for $x, y \in SHO$, we obtain
		\[\phi(x,y)\in Z_{SHO}([SHO,SHO])=Z(SHO)=0.\]
		The proof is completed.
	\end{proof}
	
\end{Lemma}

\begin{Lemma}\label{lem4.2}		
	Let $\phi$ be a skew-symmetric super-biderivation of $SHO$. For $T_H(x_ix_{i'}-x_jx_{j'})\in T_{SHO}$ and $T_H(x^{(\alpha)}x^u)\in SHO$, we have
	\[\phi(T_H(x_ix_{i'}-x_jx_{j'}),T_H(x^{(\alpha)}x^u))\in SHO_{(\alpha+\langle u\rangle)}.\]	
	
	\begin{proof}
		By lemma \ref{lem4.1}, it follows that $\phi(T_H(x_ix_{i'}-x_jx_{j'}),T_H(x_kx_{k'}-x_lx_{l'}))=0$ for any $i,j,k,l\in Y_0$ from $[T_H(x_ix_{i'}-x_jx_{j'}),$  $T_H(x_kx_{k'}-x_lx_{l'})]=0$. Note that $\d(T_H(x_ix_{i'}-x_jx_{j'}))=\bar{0}$, then for any $T_H(x^{(\alpha)}x^u)\in SHO$, it is clear that		
		\begin{align*}
			&[T_H(x_kx_{k'}-x_lx_{l'}),\phi(T_H(x_ix_{i'}-x_jx_{j'}),T_H(x^{(\alpha)}x^u))]\\
			=&(-1)^{(\d(\phi)+\bar{0})\,\bar{0}}\big(\phi(T_H(x_ix_{i'}-x_jx_{j'}),[T_H(x_kx_{k'}-x_lx_{l'}),T_H(x^{(\alpha)}x^u)])\\
			&-[\phi(T_H(x_ix_{i'}-x_jx_{j'}),T_H(x_kx_{k'}-x_lx_{l'})),T_H(x^{(\alpha)}x^u)]\big)\\
			=& (\delta_{( k'\in u)}- \alpha_k+\alpha_l-\delta_{( l'\in u)})\phi(T_H(x_ix_{i'}-x_jx_{j'}),T_H(x^{(\alpha)}x^u)).
		\end{align*}
		The proof is completed.
	\end{proof}
	
\end{Lemma}

We want to prove that all skew-symmetric super-biderivations of $SHO$ are inner. First, we need to prove the conclusion works on some elements of $SHO$. Therefore we give some specific subspaces $SHO_{(\alpha+\langle u\rangle)}$.

\begin{Lemma}\label{lem4.3}
	Let $i\in Y_0,j\in Y_1$. Then the following statements hold:
	\begin{align*}
		&(1)SHO_{(\varepsilon_i)}=\sum_{\substack{\alpha\in A,u\in B\\0\leq q <p}}\BF T_H\big((\prod_{r\in Y_0,h'\in u} x^{(\alpha_r^{\bar{q}}\varepsilon_r)}x^{(\varepsilon_h)})x^{(\varepsilon_i)}x^u\big)\cap SHO;\\
		&(2)SHO_{(\langle i' \rangle)}=\sum_{\substack{\alpha\in A,u\in B\\0\leq q <p}}\BF T_H\big((\prod_{r\in Y_0,h'\in u} x^{(\alpha_r^{\bar{q}}\varepsilon_r)}x^{(\varepsilon_h)})x_{i'}x^u\big)\cap SHO;\\
		&(3)SHO_{(\varepsilon_i+\langle j'\rangle )}=\sum_{\substack{\alpha\in A,u\in B\\0\leq q <p}}\BF T_H\big((\prod_{r\in Y_0,h'\in u} x^{(\alpha_r^{\bar{q}}\varepsilon_r)}x^{(\varepsilon_h)})x^{(\varepsilon_i)}x_{j'}x^u\big)\cap SHO;&
	\end{align*}
    where $\alpha_r^{\bar{q}}$ denotes some integer and $\alpha_r^{\bar{q}}\equiv q\pmod{p}$.
\end{Lemma}

\begin{Lemma}\label{lem4.4}
	Let $\phi$ is a skew-symmetric super-biderivation of $SHO$. For any $i,j\in Y$ and $j\neq i,i'$, there is an element $\lambda_i\in\BF$ such that
	\[\phi(T_H(x_ix_{i'}-x_jx_{j'}),T_H(x_i))=\lambda_i[T_H(x_ix_{i'}-x_jx_{j'}),T_H(x_i)],\]
	where $\lambda_i$ is dependent on $i$.
	
	\begin{proof}
		For $i'\in Y_1$, by Lemma \ref{lem4.3} (2),  we can suppose that
		\begin{align*}
			\phi(T_H(x_ix_{i'}-x_jx_{j'}),T_H(x_{i'}))
			=\sum_{\substack{\alpha\in A,u\in B\\0\leq q <p}}a(\alpha,u,i') T_H\big((\prod_{r\in Y_0,h'\in u} x^{(\alpha_r^{\bar{q}}\varepsilon_r)}x^{(\varepsilon_h)})x_{i'}x^u\big)\in SHO,
		\end{align*}				
	    where $a(\alpha,u,i')\in\BF$. For any $k\in Y\backslash\{i,i',j,j'\}$,  Lemma \ref{lem4.1} yields
	    \begin{align*}
    	  0&=(-1)^{\d(\phi)(\tau(k)+\bar{1})}\big(\phi(T_H(x_ix_{i'}-x_jx_{j'}),[T_H(x_k),T_H(x_{i'})])-[\phi(T_H(x_ix_{i'}-x_jx_{j'}),T_H(x_k)),T_H( x_{i'})]\big)\\
    	  &=[T_H(x_k),\phi(T_H(x_ix_{i'}-x_jx_{j'}),T_H(x_{i'}))]\\
    	  &=[T_H(x_k),\sum_{\substack{\alpha\in A,u\in B\\0\leq q <p}}a(\alpha,u,i') T_H\big((\prod_{r\in Y_0,h'\in u} x^{(\alpha_r^{\bar{q}}\varepsilon_r)}x^{(\varepsilon_h)})x_{i'}x^u\big)]\\
    	  &= \sum_{\substack{\alpha\in A,u\in B\\0\leq q <p}}a(\alpha,u,i') T_H\big(T_H(x_k)((\prod_{r\in Y_0,h'\in u} x^{(\alpha_r^{\bar{q}}\varepsilon_r)}x^{(\varepsilon_h)})x_{i'}x^u)\big).
	    \end{align*}
        By computing the equation, we find that $a(\alpha,u,i')=0$ if $\alpha_{k'}^{\bar{q}}>0$ or $k'\in u$. Thus we can suppose that
        \begin{align*}
            \phi(T_H(x_ix_{i'}-x_jx_{j'}),T_H(x_{i'}))
            =\sum_{\substack{\alpha\in A,\{u\}\subset \{ j'\}\\0\leq q <p}}a(\alpha,u,i') T_H\big((\prod_{h'\in u} x^{(\varepsilon_h)})x^{(\alpha_i^{\bar{q}}\varepsilon_i)}x^{(\alpha_j^{\bar{q}}\varepsilon_j)}x_{i'}x^u\big).
        \end{align*}
        Note that
        \[\phi(T_H(x_ix_{i'}-x_jx_{j'}),T_H(x_{i'}))=\phi(T_H(x_ix_{i'}-x_lx_{l'}),T_H(x_{i'})),\]
        for any $l\in Y$ and $l\neq i,i'$. From the arbitrariness of $j$, we may assume that
        \begin{align*}
        	\phi(T_H(x_ix_{i'}-x_jx_{j'}),T_H(x_{i'}))
        	=\sum_{\alpha\in A,0\leq q <p}a(\alpha,i') T_H\big(x^{(\alpha_i^{\bar{q}}\varepsilon_i)}x_{i'}\big).
        \end{align*}
        Since $\d(T_H(x_ix_{i'}-x_jx_{j'}))+\d(T_H(x_{i'}))=0$, by Lemma \ref{lem3.3}, we have
        \begin{align*}
        	0&=[[T_H(x_ix_{i'}-x_jx_{j'}),T_H(x_{i'})],\phi(T_H(x_ix_{i'}-x_jx_{j'}),T_H(x_{i'}))]\\
        	&=[T_H(x_{i'}),\phi(T_H(x_ix_{i'}-x_jx_{j'}),T_H(x_{i'}))]\\
        	&=[T_H(x_{i'}),\sum_{\alpha\in A,0\leq q <p}a(\alpha,i') T_H\big(x^{(\alpha_i^{\bar{q}}\varepsilon_i)}x_{i'}\big)]\\
        	&=\sum_{\alpha\in A,0\leq q <p}a(\alpha,i') T_H\big(T_H(x_{i'})(x^{(\alpha_i^{\bar{q}}\varepsilon_i)}x_{i'})\big).
        \end{align*}
        By computing the equation, we find that $a(\alpha,i')=0$ if $\alpha_i^{\bar{q}}>0$. Thus we can suppose that
        \begin{align*}
        	\phi(T_H(x_ix_{i'}-x_jx_{j'}),T_H(x_{i'}))
        	=a(i') T_H(x_{i'}).
        \end{align*}
       Set $\lambda_{i'}=a(i')$. Since $[T_H(x_ix_{i'}-x_jx_{j'}),T_H(x_{i'})]=T_H(x_{i'})$,
        we conclude that
        \[\phi(T_H(x_ix_{i'}-x_jx_{j'}),T_H(x_{i'}))
        =\lambda_{i'}[T_H(x_ix_{i'}-x_jx_{j'}),T_H(x_{i'})],\]
        where $\lambda_{i'}$ is dependent on $i'$.

        For $i\in Y_0$, by Lemma \ref{lem4.3} (1),  we may assume that
         \begin{align*}
         	\phi(T_H(x_ix_{i'}-x_jx_{j'}),T_H(x_i))
         	=\sum_{\substack{\alpha\in A,u\in B\\0\leq q <p}}a(\alpha,u,i) T_H\big((\prod_{r\in Y_0,h'\in u} x^{(\alpha_r^{\bar{q}}\varepsilon_r)}x^{(\varepsilon_h)})x^{(\varepsilon_i)}x^u\big)\in SHO,
         \end{align*}		
        where $a(\alpha,u,i')\in\BF$. For any $k\in Y\backslash\{i,i',j,j'\}$, Lemma \ref{lem4.1} implies that
        \begin{align*}
        	0&=(-1)^{\d(\phi)(\tau(k)+\bar{1})}\big(\phi(T_H(x_ix_{i'}-x_jx_{j'}),[T_H(x_k),T_H(x_i)])-[\phi(T_H(x_ix_{i'}-x_jx_{j'}),T_H(x_k)),T_H( x_i)]\big)\\
        	&=[T_H(x_k),\phi(T_H(x_ix_{i'}-x_jx_{j'}),T_H(x_i))]\\
        	&=[T_H(x_k),\sum_{\substack{\alpha\in A,u\in B\\0\leq q <p}}a(\alpha,u,i) T_H\big((\prod_{r\in Y_0,h'\in u} x^{(\alpha_r^{\bar{q}}\varepsilon_r)}x^{(\varepsilon_h)})x^{(\varepsilon_i)}x^u\big)]\\
        	&= \sum_{\substack{\alpha\in A,u\in B\\0\leq q <p}}a(\alpha,u,i) T_H\big(T_H(x_k)((\prod_{r\in Y_0,h'\in u} x^{(\alpha_r^{\bar{q}}\varepsilon_r)}x^{(\varepsilon_h)})x^{(\varepsilon_i)}x^u)\big).
        \end{align*}
        By computing the equation, we find that $a(\alpha,u,i')=0$ if $\alpha_{k'}^{\bar{q}}>0$ or $k'\in u$. Consequently we can suppose that
        \begin{align*}
        	\phi(T_H(x_ix_{i'}-x_jx_{j'}),T_H(x_i))
        	=\sum_{\substack{\alpha\in A,\{u\}\subset\{ i',j'\}\\0\leq q <p}}a(\alpha,u,i) T_H\big((\prod_{h'\in u} x^{(\varepsilon_h)})x^{(\alpha_i^{\bar{q}}\varepsilon_i)}x^{(\alpha_j^{\bar{q}}\varepsilon_j)}x^{(\varepsilon_i)}x^u\big).
        \end{align*}
        Similarly, from the arbitrariness of $j$, we may assume that
        \begin{align*}
        	\phi(T_H(x_ix_{i'}-x_jx_{j'}),T_H(x_i))
        	=\sum_{\substack{\alpha\in A,\{u\}\subset\{ i'\}\\0\leq q <p}}a(\alpha,u,i) T_H\big((\prod_{h'\in u} x^{(\varepsilon_h)})x^{(\alpha_i^{\bar{q}}\varepsilon_i)}x^{(\varepsilon_i)}x^u\big).
        \end{align*}
        By Lemma \ref{lem3.2}, we have
        \begin{align*}
        	0&=[\phi(T_H(x_ix_{i'}-x_jx_{j'}),T_H(x_{i'})),[T_H(x_ix_{i'}-x_jx_{j'}),T_H(x_i)]]\\
        	&=[[T_H(x_ix_{i'}-x_jx_{j'}),T_H(x_{i'})],\phi(T_H(x_ix_{i'}-x_jx_{j'}),T_H(x_i))]\\
        	&=[T_H(x_{i'}),\sum_{\substack{\alpha\in A,\{u\}\subset\{ i'\}\\0\leq q <p}}a(\alpha,u,i) T_H\big((\prod_{h'\in u} x^{(\varepsilon_h)})x^{(\alpha_i^{\bar{q}}\varepsilon_i)}x^{(\varepsilon_i)}x^u\big)]\\
        	&=\sum_{\substack{\alpha\in A,\{u\}\subset\{ i'\}\\0\leq q <p}}a(\alpha,u,i) T_H\big(T_H(x_{i'})((\prod_{h'\in u} x^{(\varepsilon_h)})x^{(\alpha_i^{\bar{q}}\varepsilon_i)}x^{(\varepsilon_i)}x^u)\big).
        \end{align*}
         By computing the equation, we find that $a(\alpha,i')=0$ if $\alpha_i^{\bar{q}}>0$ or $i'\in u$. Thus we can suppose that
        \begin{align*}
        	\phi(T_H(x_ix_{i'}-x_jx_{j'}),T_H(x_i))
        	=a(i') T_H(x_i).
        \end{align*}
       Set $\lambda_{i}=-a(i)$. Since $[T_H(x_ix_{i'}-x_jx_{j'}),T_H(x_i)]=-T_H(x_i)$,
       we obtain 
        \[\phi(T_H(x_ix_{i'}-x_jx_{j'}),T_H(x_i))
        =\lambda_i[T_H(x_ix_{i'}-x_jx_{j'}),T_H(x_i)],\]
        where $\lambda_i$ is dependent on $i$.
	\end{proof}

\end{Lemma}

\begin{Lemma}\label{lem4.5}
	All $\BZ_2$-homogeneous skew-symmetric super-biderivations of $SHO$ are even.
	
	\begin{proof}
			Suppose $\phi$ is a $\BZ_2$-homogeneous skew-symmetric super-biderivation of $SHO$. Due to Lemma \ref{lem4.4}, we know that  $\phi(T_H(x_ix_{i'}-x_jx_{j'}),T_H(x_i))$ and $[T_H(x_ix_{i'}-x_jx_{j'}),T_H(x_i)]$ have the same $\BZ_2$-degree, so $\phi$ and $[\cdot,\cdot]$ have the same $\BZ_2$-degree. Then $\phi$ is even.
		\end{proof}

\end{Lemma}

\begin{Corollary}\label{cor4.6}
	All skew-symmetric super-biderivations of $SHO$ are $\BZ_2$-homogeneous, and further even.
\end{Corollary}

\begin{Lemma}\label{lem4.7}
	Let $\phi$ be a skew-symmetric super-biderivation of $SHO$. For any $i,j\in Y_0$ and $i\neq j$, there is an element $\lambda\in\BF$ such that
	\[\phi(T_H(x_ix_{i'}-x_jx_{j'}),T_H(x_ix_{j'}))=\lambda[T_H(x_ix_{i'}-x_jx_{j'}),T_H(x_ix_{j'})],\]
	where $\lambda$ depends on neither $i$ nor $j$.
	
	\begin{proof}
		By Lemma \ref{lem4.3} (3),  we can suppose that
		\begin{align*}
			\phi(T_H(x_ix_{i'}-x_jx_{j'}),T_H(x_ix_{j'}))
			=\sum_{\substack{\alpha\in A,u\in B\\0\leq q <p}}a(\alpha,u,i,j') T_H\big((\prod_{r\in Y_0,h'\in u} x^{(\alpha_r^{\bar{q}}\varepsilon_r)}x^{(\varepsilon_h)})x^{(\varepsilon_i)}x_{j'}x^u\big)\in SHO,
		\end{align*}				
		where $a(\alpha,u,i,j')\in\BF$. According to Lemma \ref{lem4.1}, for any $k\in Y\backslash\{i,i',j,j'\}$, we have
		\begin{align*}
			0&=\phi(T_H(x_ix_{i'}-x_jx_{j'}),[T_H(x_k),T_H(x_ix_{j'})])-[\phi(T_H(x_ix_{i'}-x_jx_{j'}),T_H(x_k)),T_H(x_ix_{j'})]\\
			&=[T_H(x_k),\phi(T_H(x_ix_{i'}-x_jx_{j'}),T_H(x_ix_{j'}))]\\
			&=[T_H(x_k),\sum_{\substack{\alpha\in A,u\in B\\0\leq q <p}}a(\alpha,u,i,j') T_H\big((\prod_{r\in Y_0,h'\in u} x^{(\alpha_r^{\bar{q}}\varepsilon_r)}x^{(\varepsilon_h)})x^{(\varepsilon_i)}x_{j'}x^u\big)]\\
			&= \sum_{\substack{\alpha\in A,u\in B\\0\leq q <p}}a(\alpha,u,i,j') T_H\big(T_H(x_k)((\prod_{r\in Y_0,h'\in u} x^{(\alpha_r^{\bar{q}}\varepsilon_r)}x^{(\varepsilon_h)})x^{(\varepsilon_i)}x_{j'}x^u)\big).
		\end{align*}
		By computing the equation, we find that $a(\alpha,u,i,j')=0$ if $\alpha_{k'}^{\bar{q}}>0$ or $k'\in u$. Thus we can suppose that
		\begin{align*}
			\phi(T_H(x_ix_{i'}-x_jx_{j'}),T_H(x_ix_{j'}))
			=\sum_{\substack{\alpha\in A,\{u\}\in\{i'\}\\0\leq q <p}}a(\alpha,u,i,j') T_H\big((\prod_{h'\in u} x^{(\varepsilon_h)})x^{(\alpha_i^{\bar{q}}\varepsilon_i)}x^{(\alpha_j^{\bar{q}}\varepsilon_j)}x^{(\varepsilon_i)}x_{j'}x^u\big).
		\end{align*}
	    By Lemmas \ref{lem3.2} and  \ref{lem4.4}, we have
	    \begin{align*}
	    	2\lambda_{i'}T_H(x_{j'})&=[\lambda_{i'}[T_H(x_ix_{i'}-x_jx_{j'}),T_H(x_{i'})],[T_H(x_ix_{i'}-x_jx_{j'}),T_H(x_ix_{j'})]]\\
	    	&=[\phi(T_H(x_ix_{i'}-x_jx_{j'}),T_H(x_{i'})),[T_H(x_ix_{i'}-x_jx_{j'}),T_H(x_ix_{j'})]]\\
	    	&=[[T_H(x_ix_{i'}-x_jx_{j'}),T_H(x_{i'})],\phi(T_H(x_ix_{i'}-x_jx_{j'}),T_H(x_ix_{j'}))]\\
	    	&=[T_H(x_{i'}),\sum_{\substack{\alpha\in A,\{u\}\in\{i'\}\\0\leq q <p}}a(\alpha,u,i,j') T_H\big((\prod_{h'\in u} x^{(\varepsilon_h)})x^{(\alpha_i^{\bar{q}}\varepsilon_i)}x^{(\alpha_j^{\bar{q}}\varepsilon_j)}x^{(\varepsilon_i)}x_{j'}x^u\big)]\\
	    	&= \sum_{\substack{\alpha\in A,\{u\}\in\{i'\}\\0\leq q <p}}a(\alpha,u,i,j') T_H\big(T_H(x_{i'})((\prod_{h'\in u} x^{(\varepsilon_h)})x^{(\alpha_i^{\bar{q}}\varepsilon_i)}x^{(\alpha_j^{\bar{q}}\varepsilon_j)}x^{(\varepsilon_i)}x_{j'}x^u)\big).
	    \end{align*}		
		By computing the equation, we find that $a(\alpha,u,i,j')=0$ if $\alpha_i^{\bar{q}}>0$, $\alpha_j^{\bar{q}}>0$ or $i'\in u$. Thus we may assume that
		\begin{align*}
			\phi(T_H(x_ix_{i'}-x_jx_{j'}),T_H(x_ix_{j'}))
			=a(i,j')T_H(x_ix_{j'})=-2\lambda_{i'}T_H(x_ix_{j'}).
		\end{align*}
		Since $[T_H(x_ix_{i'}-x_jx_{j'}),T_H(x_ix_{j'})]=-2T_H(x_ix_{j'})$, we conclude that
		\[\phi(T_H(x_ix_{i'}-x_jx_{j'}),T_H(x_ix_{j'}))
		=\lambda_{i'}[T_H(x_ix_{i'}-x_jx_{j'}),T_H(x_ix_{j'})],\]
		 Lemmas \ref{lem3.2} and  \ref{lem4.4}  imply that 
		\begin{align*}			
			0=&[\phi(T_H(x_jx_{j'}-x_ix_{i'}),T_H(x_j)),[T_H(x_ix_{i'}-x_jx_{j'}),T_H(x_ix_{j'})]]\\
			    &-[[T_H(x_jx_{j'}-x_ix_{i'}),T_H(x_j)],\phi(T_H(x_ix_{i'}-x_jx_{j'}),T_H(x_ix_{j'}))]\\			
			=&(\lambda_j-\lambda_{i'})[T_H(x_j),2T_H(x_ix_{j'})]\\
			=&2(\lambda_j-\lambda_{i'})T_H(x_i).
		\end{align*}		
		Since $T_H(x_i)\neq 0$, we obtain $\lambda_j=\lambda_{i'}$. From the arbitrariness of $i'$ and $j$, we can claim that $\lambda_1=\lambda_2=\dots=\lambda_s$. Set $\lambda:=\lambda_1=\dots=\lambda_s$. Then we can conclude that for any $i,j\in Y_0$ and $i\neq j$, there is an element $\lambda\in\BF$ such that
		\[\phi(T_H(x_ix_{i'}-x_jx_{j'}),T_H(x_ix_{j'}))=\lambda[T_H(x_ix_{i'}-x_jx_{j'}),T_H(x_ix_{j'})],\]
		where $\lambda$ depends on neither $i$ nor $j$.
	\end{proof}
	
\end{Lemma}

Then the conclusions of Lemmas \ref{lem4.4} and \ref{lem4.7} can be rewritten as
\begin{align*}
	\phi(T_H(x_ix_{i'}-x_jx_{j'}),T_H(x_i))&=\lambda[T_H(x_ix_{i'}-x_jx_{j'}),T_H(x_i)],\\
	\phi(T_H(x_ix_{i'}-x_jx_{j'}),T_H(x_{i'}))&=\lambda[T_H(x_ix_{i'}-x_jx_{j'}),T_H(x_{i'})],\\
	\phi(T_H(x_ix_{i'}-x_jx_{j'}),T_H(x_ix_{j'}))&=\lambda[T_H(x_ix_{i'}-x_jx_{j'}),T_H(x_ix_{j'})],
\end{align*}
where $i,j\in Y_0$, and $i\neq j$.

\begin{Theorem}
	Let $SHO$ be the special odd Hamiltonian superalgebra $SHO(n,n;\underline{t})$ over the basic field $\BF$ of characteristic $p>2$. Then all skew-symmetric super-biderivations of $SHO$ are inner, i.e.
	\[\BDer(SHO)=\IBDer(SHO).\]
	\begin{proof}
		Suppose $\phi$ is a skew-symmetric super-biderivation of $SHO$. For any $T_H(f),T_H(g)\in SHO$, by virtue of   Lemmas \ref{lem3.2} and  \ref{lem4.4}, for any $i\in Y$ we find that
		\begin{align*}
			0=&[\phi(T_H(x_ix_{i'}-x_jx_{j'}),T_H(x_i)),[T_H(f),T_H(g)]]\\
			    &-[[T_H(x_ix_{i'}-x_jx_{j'}),T_H(x_i)],\phi(T_H(f),T_H(g))]\\
			 =&[\lambda[T_H(x_ix_{i'}-x_jx_{j'}),T_H(x_i)],[T_H(f),T_H(g)]]\\
			 &-[[T_H(x_ix_{i'}-x_jx_{j'}),T_H(x_i)],\phi(T_H(f),T_H(g))]\\
			 =&-[T_H(x_i),\lambda[T_H(f),T_H(g)]-\phi(T_H(f),T_H(g))]\\
			 =&-[D_{i'},\lambda[T_H(f),T_H(g)]-\phi(T_H(f),T_H(g))].
		\end{align*}
	    Since $Z_{SHO}(SHO_{-1})=\{E\in SHO\mid [E,D_i]=0,\forall\ i\in Y\}=SHO_{-1}$, we obtain
	        \[\phi(T_H(f),T_H(g))-\lambda[T_H(f),T_H(g)]=\sum_{l\in Y}b_lD_l,\]
	    where $b_l\in \BF$.   Lemmas \ref{lem3.2} and   \ref{lem4.7} yield
	    \begin{align*}
	    	0=&[\phi(T_H(x_ix_{i'}-x_jx_{j'}),T_H(x_ix_{j'})),[T_H(f),T_H(g)]]\\
	    	&-[[T_H(x_ix_{i'}-x_jx_{j'}),T_H(x_ix_{j'})],\phi(T_H(f),T_H(g))]\\
	    	=&[\lambda[T_H(x_ix_{i'}-x_jx_{j'}),T_H(x_ix_{j'})],[T_H(f),T_H(g)]]\\
	    	&-[[T_H(x_ix_{i'}-x_jx_{j'}),T_H(x_ix_{j'})],\phi(T_H(f),T_H(g))]\\
	    	=&-2[T_H(x_ix_{j'}),\lambda[T_H(f),T_H(g)]-\phi(T_H(f),T_H(g))]\\
	    	=&2[T_H(x_ix_{j'}),\sum_{l\in Y}b_lD_l]\\
	    	=&-2(b_{j'}T_H(x_i)+b_iT_H(x_{j'})).
	    \end{align*}
         As $T_H(x_i)\neq0$ and $T_H(x_{j'})\neq0$,   $b_i=b_{j'}=0$. Then $b_l=0,\ \forall\ l\in Y$. Hence, for any $T_H(f),T_H(g)\in SHO$, we  conclude that
             \[\phi(T_H(f),T_H(g))=\lambda[T_H(f),T_H(g)].\]
         Thus $\phi$ is an inner super-biderivation. The assertion follows from the arbitrariness of $\phi$.
	\end{proof}
\end{Theorem}

%

{\bf Competing interests} The authors have no relevant financial or non-financial interests to disclose.

{\bf Data availability} All data generated or analysed during this study are included in this published
article.

{\bf Author contributions} Xiaoning Xu contributed to the study conception and design. The first
draft of the manuscript was written by Da Xu and all authors commented on previous versions
of the manuscript. All authors read and approved the final manuscript.


\begin{thebibliography}{26}
	\setlength{\itemsep}{0mm} \small
	
	\bibitem{a1} Y. Chang, L.Y. Chen, Biderivations and linear commuting maps on the restricted Cartan-type Lie algebras $W(n;\underline{1})$ and $S(n;\underline{1})$. \emph{Linear Multilinear Algebra}, {\bf 67}(8) (2019) 1625-1636.
	
	\bibitem{a2} Y. Chang, L.Y. Chen, X. Zhou, Biderivations and linear commuting maps on the restricted Cartan-type Lie algebras $H(n;\underline{1})$. \emph{Comm. Algebra}, {\bf 47}(3) (2019) 1311-1326.
	
	\bibitem{a3} Y. Chang, L.Y. Chen, Y. Cao, Super-biderivations of the generalized Witt Lie superalgebra $W(m, n;\underline{t})$. \emph{Linear Multilinear Algebra}, {\bf 69}(2) (2021) 233-244.
	
	\bibitem{a4} Z.X. Chen, Biderivations and linear commuting maps on simple generalized Witt algebras over a field. \emph{Electron. J. Linear Algebra}, {\bf 31} (2016) 1-12.
	
	\bibitem{a5} J.Y. Fu, Q.C. Zhang, C.P. Jiang, The Cartan-type modular Lie superalgebra KO. \emph{Comm. Algebra}, {\bf34}(1) (2006) 107-128.
	
	\bibitem{a6} G.Z. Fan, X.S. Dai, Super-biderivations of Lie superalgebras. \emph{Linear Multilinear algebra}, {\bf 65}(1) (2017) 58-66.
	
	\bibitem{a7}B.L. Guan, L.Y. Chen, Derivations of the even part of contact Lie superalgebra. \emph{J. Pure Appl. Algebra},  {\bf 216} (2012) 1454--1466.
	
	\bibitem{a8} X. Han, D.Y. Wang, C.G. Xia, Linear commuting maps and biderivations on the Lie algebras $W (a, b)$. \emph{J. Lie theory}, {\bf 26}(3) (2016) 777-786.
	
	\bibitem{a9} V. G. Kac, Lie superalgebras. \emph{Adv. Math.}, {\bf 26} (1977) 8-96.
	
	\bibitem{a10} V.G. Kac, Classification of infinite-dimensional simple linearly compact Lie superalgebras. \emph{Adv. Math.}, {\bf 139} (1998) 1-55.
	
	\bibitem{a11} D. Leites, Towards classification of simple finite dimensional modular Lie superalgebras. \emph{J. Prime Res. Math.}, {\bf 3} (2007) 101-110.
	
	\bibitem{a12} W.D. Liu, Y.H. He, Finite-dimensional special odd Hamiltonian superalgebras in prime characteristic. \emph{Commun. Contemp. Math.}, {\bf11}(04) (2009) 523-546.
	
	\bibitem{a13} W.D. Liu, Y.Z. Zhang, X.L. Wang, The derivation algebra of the Cartan-type Lie superalgebra $HO$. \emph{J. Algebra}, {\bf 273}(1) (2004) 176-205.
	
	\bibitem{a14} Q. Mu, Y.Z. Zhang, Infinite-dimensional modular special odd contact superalgebras. \emph{J. Pure Appl. Algebra}, {\bf214}(8) (2010) 1456-1468.
	
	\bibitem{a15}B. Shu, C.W. Zhang, Restricted representations of the Witt superalgebras. \emph{J. Algebra}, {\bf 324}(4) (2010) 652-672.
	
	\bibitem{a16} H.Z. Sun, Q.Z. Han, A survey of Lie superalgebras. \emph{Adv. Phys.(PRC)}, {\bf 1} (1983) 81-125 (in Chinses).
	
	\bibitem{a17} H.Z. Sun, Q.Z. Han, Lie Algebras, Lie Superalgebras and their Applications in Physics. \emph{Peking Univ.
		Press, Beijing}, (1999) (in Chinese).
	
	\bibitem{a18} M. Scheunert, Theory of Lie superalgebras. \emph{Lecture Notes in Math.}, {\bf 716} (1979) 3001-454.
	
	\bibitem{a19} X.M. Tang, Biderivations of finite-dimensional complex simple Lie algebras. \emph{Linear Multilinear algebra}, {\bf 66}(2) (2018) 250-259.
	
	\bibitem{a20} D.Y. Wang, X.X. Yu, Z.X. Chen, Biderivations of the parabolic subalgebras of simple Lie algebras. \emph{Comm. Algebra}, {\bf 39}(11) (2011) 4097-4104.
	
	\bibitem{a21} D.Y. Wang, X.X. Yu, Biderivations and linear commuting maps on the Schr$\ddot{\text{o}}$dinger-Virasoro Lie algebra. \emph{Comm. Algebra}, {\bf 41}(6) (2013) 2166-2173.
	
	\bibitem{a22} C.G. Xia, D.Y. Wang, X. Han, Linear super-commuting maps and super-biderivations on the super-Virasoro algebras. \emph{Comm. Algebra}, {\bf 44}(12) (2016) 5342-5350.
	
	\bibitem{a23} J.X. Yuan, X.M. Tang, Super-biderivations of classical simple Lie superalgebras. \emph{Aequationes Math.}, {\bf 92}(1) (2018) 91-109.
	
	\bibitem{a24} X.D. Zhao, Y. Chang, X. Zhou, L.Y. Chen, Super-biderivations of the contact Lie superalgebra $K(m, n;\underline{t})$. \emph{Comm. Algebra}, {\bf 488} (2020) 3237-3248.
	
	\bibitem{a25} Y.Z. Zhang, Finite-dimensional Lie superalgebras of Cartan type over fields of prime characteristic. \emph{Chinese Sci. Bull}, {\bf 42}(9) (1997) 720-724.
	
\end{thebibliography}
\end{document}